\documentclass[11pt]{article}
\usepackage{graphicx, colordvi}
\usepackage{dsfont}
\usepackage{epsfig}
\usepackage{mathrsfs}
\usepackage{amssymb,amsfonts,amsmath,amsthm,cite,color}
\usepackage{dsfont}
\usepackage{CJKutf8}
\usepackage{epsfig}
\usepackage{mathrsfs}
\usepackage{longtable}
\usepackage{enumerate}
\usepackage{geometry}
\usepackage{hyperref}
\hypersetup{
	colorlinks=true,
	linkcolor=cyan,
	filecolor=blue,
	urlcolor=red,
	citecolor=green,
}

\newtheorem{theo}{Theorem}

\newtheorem{lem}[theo]{Lemma}

\newtheorem{coj}[theo]{Conjecture}
\makeatletter \@addtoreset{equation}{section}
\@addtoreset{theo}{section} \makeatother

\newcommand{\bN} { {\mathbb{N}}}

\newcommand{\bZ} { {\mathbb{Z}}}

\def\qed{\hfill \rule{4pt}{7pt}}
\def\pf{\noindent {\it Proof.} }

\author{}
\title{}

\begin{document}
\begin{center}
 {\large \bfseries Arithmetic properties of generalized Delannoy polynomials and Schr\"oder polynomials}
\end{center}

\begin{center}
 {Lin-Yue Li}$^{1}$ and {Rong-Hua Wang}$^{2}$

   $^{1,2}$School of Mathematical Sciences\\
   Tiangong University \\
   Tianjin 300387, P.R. China\\
   $^{1}$2330141409@tiangong.edu.cn\\
   $^{2}$wangronghua@tiangong.edu.cn \\[10pt]
\end{center}

\vskip 6mm \noindent {\bf Abstract.}
Let $n$ be any nonnegative integer and
\[
D_n^{(h)}(x)=\sum_{k=0}^{n}\binom{n+k}{2k}^{h}\binom{2k}{k}^{h}{x}^{k}
\text{ and }
S_{n}^{(h)}(x)=\sum_{k=0}^{n}\binom{n+k}{2k}^{h}C_{k}^{h}{x}^{k}
\]
be the generalized Delannoy polynomials and Schr\"oder polynomials respectively.
Here $C_k$ is the Catalan number and $h$ is a positive integer.
In this paper, we prove that
\begin{align*}
& \frac{(2,n)}{n(n+1)(n+2)}
\sum_{k=1}^{n}k^a(k+1)^a(2k+1)D_{k}^{(h)}(x)^{m}\in\bZ[x],\\
&\frac{(2,hm-1,n)}{n(n+1)(n+2)}
\sum_{k=1}^{n}(-1)^{k}k^a(k+1)^a(2k+1)D_{k}^{(h)}(x)^{m}\in\bZ[x],\\
&\frac{(2,n)}{n(n+1)(n+2)}
\sum_{k=1}^{n}k^a(k+1)^a(2k+1)S_{k}^{(h)}(x)^{m}\in\bZ[x],\\
&\frac{(2,m-1,n)}{n(n+1)(n+2)}
\sum_{k=1}^{n}(-1)^{k}k^a(k+1)^a(2k+1)S_{k}^{(h)}(x)^{m}\in\bZ[x].
\end{align*}
Taking $a=1$ will confirm some of Z.-W. Sun's conjectures.

	\section{Introduction}
Let $\bN$ and $\bZ^{+}$ denote the set of nonnegative and positive integers respectively. Let $n\in\bN$, the $n$-th central Delannoy numbers $D_n$ and the $n$-th large Schr\"oder numbers $S_n$ \cite{RP1999} are defined respectively by
\begin{equation*}
D_{n}=\sum\limits_{k=0}^{n}\binom{n}{k}\binom{n+k}{k}
     =\sum_{k=0}^{n}\binom{n+k}{2k}\binom{2k}{k}
\end{equation*}
and
\begin{equation*}
S_n=\sum_{k=0}^{n}\binom{n}{k}\binom{n+k}{k}\frac{1}{k+1}
   =\sum_{k=0}^{n}\binom{n+k}{2k}C_k,
\end{equation*}
where $C_k$ denotes the Catalan number $\binom{2k}{k}/(k+1)$.

There are many combinatorial objects that are counted by $D_n$ and $S_n$.
For example, $D_n$ is the number of lattice paths from the point $(0, 0)$ to $(n, n)$ with steps choosing from the set $\{(1, 0), (0, 1), (1, 1)\}$, while $S_n$ counts lattice paths satisfying the same constraints together with the limitation that never rise above the line $y=x$.
One can consult items A001850 and A006318 in the OEIS \cite{The OEIS Foundation Inc}
for other classical combinatorial interpretations.

Sun \cite{Sun2011,Sun2014} discovered many amazing arithmetic properties of $D_n$ and $S_n$.
For example he proved
\[
\sum_{k=1}^{p-1}\frac{D_k}{k^2}\equiv (-1)^{(p-1)/2}2E_{p-3}
\pmod p \text{\quad and \quad}
\sum_{k=1}^{p-1}\frac{S_k}{6^k}\equiv 0\pmod p
\]
for any prime $p>3$, where $E_k$ are the Euler numbers,
and conjectured that
\[
\sum_{k=1}^{p-1}D_kS_k\equiv -2p\sum_{k=1}^{p-1}\frac{3+(-1)^k}{k} \pmod {p^4}
\]
which was confirmed by Liu in \cite{Liu2016}.
Sun \cite{Sun2014} also showed that
\[
\frac{1}{n^2}\sum_{k=0}^{n-1}(2k+1)D_k^2\in\bZ
\text{ for } n\in\bZ^{+}.
\]
In 2012, Guo and Zeng\cite{GuoZeng2012a} derived
\begin{equation*}
	\sum_{k=0}^{n-1}\varepsilon^k(2k+1)k^a(k+1)^aD_{k}
\equiv\sum_{k=0}^{n-1}\varepsilon^k(2k+1)^{2a+1}D_{k}\equiv 0 \pmod {n},
\end{equation*}
where $n\geq1$, $a\geq0$, and $\varepsilon=\pm 1$.
In 2017, Cao and Pan\cite{CaoPan2017} showed that,
\begin{equation*}
	S_{n+2^{\alpha}}\equiv S_n+2^{\alpha+1}\pmod{2^{\alpha+2}}
\end{equation*}
for $n\geq1$ and $\alpha\geq1$.

Motivated by the deﬁnition of Delannoy numbers and Schr\"oder numbers, Sun \cite{ Sun2012a, Sun2012b} defined the Delannoy polynomials $D_{n}(x)$ and the large Schr\"oder
polynomials $S_n(x)$ as
\begin{equation*}
	D_{n}(x)=\sum\limits_{k=0}^{n}\binom{n}{k}\binom{n+k}{k}x^k\quad \text{and} \quad 	S_n(x)=\sum_{k=0}^{n}\binom{n}{k}\binom{n+k}{k}\frac{x^k}{k+1}.
\end{equation*}
Clearly, $D_{n}(1)$ and $S_n(1)$ reduce to the $n$-th central Delannoy numbers and the $n$-th large Schr\"oder numbers respectively.
The Delannoy polynomials and large Schr\"oder polynomials also have nice arithmetic properties.
For example, Sun\cite{Sun2014} proved that
\begin{equation*}
	\sum_{k=0}^{p-1}D_k(x)^2\equiv\left( \frac{x(x+1)}{p}\right)\pmod p
\end{equation*}
for any odd prime $p$ and integer $x$, where $\left(\frac{\bullet}{p}\right)$ denotes the Legendre symbol.
Guo\cite{Guo2015} confirmed Sun's conjecture in \cite{Sun2014} that for any prime $p$ and integer $x$ satisfying $p\nmid x(x+1)$,
\begin{equation*}
	\sum_{k=0}^{p-1}(2k+1)D_k(x)^3\equiv p\left( \frac{-4x-3}{p}\right) \pmod  {p^2}
\end{equation*}
and
\begin{equation*}
     \sum_{k=0}^{p-1}(2k+1)D_k(x)^4\equiv p\pmod  {p^2}.
\end{equation*}
Liu\cite{Liu2016} proved that
\begin{equation*}
	\sum_{k=1}^{p-1}D_k(x)S_k(x)\equiv0\pmod p,
\end{equation*}
where $p$ is an odd prime and $x$ is an integer not divisible by $p$.
In 2018, Sun \cite{Sun2018} found
\begin{align*}
	\frac{1}{n}\sum_{k=0}^{n-1}D_{k}(x)s_{k+1}(x)\in\bZ[x(x+1)] \quad\text{for all $n\in\bZ^{+}$}
\end{align*}
and
\begin{align*}
	\sum_{k=0}^{p-1}D_{k}(x)s_{k+1}(x)\equiv0\pmod{p^{2}}
\end{align*}
for any odd prime $p$ and $p$-adic integer $x\not\equiv0,-1\pmod p$, where $s_n(x)$ is the little Schr\"oder
polynomials given by $s_n(x)=\sum_{k=1}^{n}\frac{1}{n}\binom{n}{k}\binom{n}{k-1}x^{k-1}(x+1)^{n-k}$.
In 2024, Jia, Wang and Zhong\cite{JiaWang 2024b} showed that
 \begin{equation}\label{eq:JWZ}
 	\sum_{k=0}^{p-1}(2k+1)^{2a+1}\varepsilon^{k}S_k(x)\equiv1 \pmod p
 \end{equation}
for any odd prime $p$, nonnegative integers $a\in\bN$, $\varepsilon\in\{-1,1\}$ and integer $x$ with $p\nmid x(x+1)$.

The polynomials $D_n(x)$ and $S_n(x)$ can be further generalized.
In 2012, Guo and Zeng \cite{GuoZeng2012b} introduced the Schmidt polynomials
\begin{equation}\label{eq:Dn(x)}
D_n^{(h)}(x)=\sum_{k=0}^{n}\binom{n}{k}^{h}\binom{n+k}{k}^{h}x^{k}
            =\sum_{k=0}^{n}\binom{n+k}{2k}^{h}\binom{2k}{k}^{h}{x}^{k}.
\end{equation}
Pan \cite{Pan2014} proved that, for any positive integers $h$, $m$ and $n$, there holds
\begin{equation*}
	\sum_{k=0}^{n-1}\varepsilon^{k}(2k+1)D_{k}^{(h)}(x)^{m}\equiv 0 \pmod n,
\end{equation*}
where $\varepsilon=\pm 1$.
In 2018, Chen and Guo \cite{chenguo2018} defined the multi-variable Schmidt polynomials as
\begin{equation*} S_{n}^{(h)}(x_{0},x_1,\ldots,x_{n})
=\sum_{k=0}^{n}\binom{n+k}{2k}^{h}\binom{2k}{k}x_{k}, \quad h\in\bZ^{+}
	\end{equation*}
and proved that for $a\in\bN$ and $h,m,n\in\bZ^{+}$, one have
	\begin{equation}\label{eq:CHenGuok^a(k+1)^a(2k+1)}
\sum_{k=0}^{n-1}\varepsilon^{k}k^a(k+1)^a(2k+1)S^{(h)}_k(x_0,x_1,\ldots,x_n)^m
\equiv\ 0 \pmod{n},
	\end{equation}
where $\varepsilon=\pm 1$.

In 2022, Sun\cite{Sun2022} defined
	\begin{equation}\label{eq:S_{k}^{(h)}(x)}
	S_{n}^{(h)}(x)=\sum_{k=0}^{k}\binom{n+k}{2k}^{h}C_{k}^{h}{x}^{k}, n\in\bN
	\end{equation}
for $h\in\bZ^{+}$
and provided many interesting conjectures involving $S_k^{(h)}(x)$ and $D_k^{(h)}(x)$.
\begin{coj}[Conjecture 5.3 of Sun\cite{Sun2022}]\label{th:2022sun Conjecture 5.3(2)}
	For any ${h,m,n}\in\bZ^{+}$,
	
	\textnormal{(i)} we have
	\begin{equation*}
		\frac{(2,n)}{n(n+1)(n+2)}\sum_{k=1}^{n}k(k+1)(2k+1)D_{k}^{(h)}(x)^{m}\in\bZ[x]
	\end{equation*}
	and
	\begin{equation*}
		\frac{(2,hm-1,n)}{n(n+1)(n+2)}
		\sum_{k=1}^{n}(-1)^{k}k(k+1)(2k+1)D_{k}^{(h)}(x)^{m}\in\bZ[x].
	\end{equation*}
\textnormal{(ii)} we have
\begin{equation*}
	\frac{(2,n)}{n(n+1)(n+2)}\sum_{k=1}^{n}k(k+1)(2k+1)S_{k}^{(h)}(x)^{m}\in\bZ[x]
\end{equation*}
and
\begin{equation*}	\frac{(2,m-1,n)}{n(n+1)(n+2)}\sum_{k=1}^{n}(-1)^{k}k(k+1)(2k+1)S_{k}^{(h)}(x)^{m}\in\bZ[x].
\end{equation*}
\end{coj}

In this paper, we confirm and generalize the above four conjectures.
More precisely, we proved that
\begin{align*}
\frac{(2,n)}{n(n+1)(n+2)} & \sum_{k=1}^{n}k^a(k+1)^a(2k+1)D_{k}^{(h)}(x)^{m}\in\bZ[x],\\
\frac{1}{n(n+1)(n+2)}     &
                            \sum_{k=1}^{n}(-1)^{k}k^a(k+1)^a(2k+1)D_{k}^{(h)}(x)^{m}\in\bZ[x]
			                    \ for\ h>1, \\
\frac{(2,m-1,n)}{n(n+1)(n+2)}&\sum_{k=1}^{n}(-1)^{k}k^a(k+1)^a(2k+1)D_{k}(x)^{m}\in\bZ[x],\\
\frac{(2,n)}{n(n+1)(n+2)}    & \sum_{k=1}^{n}k^a(k+1)^a(2k+1)S_{k}^{(h)}(x)^{m}\in\bZ[x], \\
\frac{(2,m-1,n)}{n(n+1)(n+2)} &
                              \sum_{k=1}^{n}(-1)^{k}k^a(k+1)^a(2k+1)S_{k}^{(h)}(x)^{m}\in\bZ[x]
\end{align*}
for ${h,m,n,a}\in\bZ^{+}$.

The rest of the paper is organized as follows.
In section 2, we provide proofs of the first three results.
The remained two will be discussed in section 3.

\section{Congruences involving $D_{k}^{(h)}(x)^{m}$}
In this section, we study the arithmetic properties of $D_{k}^{(h)}(x)^{m}$.
\begin{theo}\label{TM:Dn1}
Let $D_k^{(h)}(x)$ be given as in \eqref{eq:Dn(x)} and $\varepsilon\in\{-1,1\}$.
For any  ${h,m,n,a}\in\bZ^{+}$, we have
\begin{equation}\label{eq: (2,n){n(n+1)(n+2)}}
\frac{(2,n)}{n(n+1)(n+2)}
\sum_{k=1}^{n}\varepsilon^{k}k^a(k+1)^a(2k+1)D_{k}^{(h)}(x)^{m}\in\bZ[x].
\end{equation}

\end{theo}
\pf
Let $x_k=\binom{2k}{k}^{h-1}x^k$ in \eqref{eq:CHenGuok^a(k+1)^a(2k+1)} gives
\begin{equation}\label{eq:CGk^a(k+1)^a(2k+1)}
n\mid\sum_{k=0}^{n-1}\varepsilon^{k}k^a(k+1)^a(2k+1)D_{k}^{(h)}(x)^{m}.
\end{equation}
Note that $a\geq 1$, this leads to
\begin{equation}\label{eq:n}
	n\mid \sum_{k=1}^{n}\varepsilon^{k}k^a(k+1)^a(2k+1)D_{k}^{(h)}(x)^{m}
\end{equation}
since the $n$-th term in the summation contains $n^a$.
Let $n$ be $n+1$ and $n+2$ in \eqref{eq:n}, similar discussions leads to
\begin{equation}\label{eq:n+1}
	(n+1)\mid \sum_{k=1}^{n}\varepsilon^{k}k^a(k+1)^a(2k+1)D_{k}^{(h)}(x)^{m}
\end{equation}
and
\begin{equation}\label{eq:n+2}
	(n+2)\mid \sum_{k=1}^{n}\varepsilon^{k}k^a(k+1)^a(2k+1)D_{k}^{(h)}(x)^{m}.
\end{equation}
Then  \eqref{eq: (2,n){n(n+1)(n+2)}} follows from \eqref{eq:n}--\eqref{eq:n+2} by noting that
$[n(n+1),(n+2)]=\frac{n(n+1)(n+2)}{(2,n)}.$
\qed

When $\varepsilon=-1$, a stronger version can be given.
\begin{theo}\label{TM:Dn2}
For any  ${h,m,n,a}\in\bZ^{+}$, we have
\begin{equation}\label{eq:1}	\frac{1}{n(n+1)(n+2)}\sum_{k=1}^{n}(-1)^{k}k^a(k+1)^a(2k+1)D_{k}^{(h)}(x)^{m}\in\bZ[x]
\ for\ h>1,
\end{equation}
and
\begin{equation}\label{eq:3Dk}
\frac{(2,m-1,n)}{n(n+1)(n+2)}\sum_{k=1}^{n}(-1)^{k}k^a(k+1)^a(2k+1)D_{k}(x)^{m}\in\bZ[x].
\end{equation}
\end{theo}

The proofs of \eqref{eq:1} and \eqref{eq:3Dk} are much more complicated.
To proceed, some preliminaries are needed.	
\begin{lem}\label{lem: D_{k}^{(h)}(x)}
For $n,u,\ell\in\bN$, we have
\begin{equation}\label{eq:2k+1}
\sum_{k=0}^{n}(-1)^{k}(2k+1)\binom{k+\ell+u}{2\ell+2u}
=(-1)^{n}(n+1+\ell+u)\binom{n+\ell+u}{2\ell+2u}.
\end{equation}
\end{lem}
\noindent\emph{Proof.}
By Gosper’s algorithm \cite{Gosper1978}, we have
\begin{equation}\label{eq:(k-l)(k^{2}-l-1)}
(-1)^{k}(2k+1)\binom{k+\ell+u}{2\ell+2u}=\Delta_{k}\left((-1)^{k+1}(k-\ell-u)\binom{k+\ell+u}{2\ell+2u}\right),
\end{equation}
where $\Delta_k$ is the difference operator
(that is, $\Delta_k F(k) = F(k + 1)-F(k)$).
Summing both sides of equality \eqref{eq:(k-l)(k^{2}-l-1)} over
$k$ from $0$ to $n$ leads to the equality \eqref{eq:2k+1}.
\qed
\begin{lem}\label{lem:{2}{n+2}}
		For any $n,\ell\in\bZ^{+}$, we have
	\begin{equation}\label{eq:l}
			\frac{2}{n+2}\binom{n-1}{\ell-1}\binom{n+\ell+1}{\ell}\in\bZ.
	\end{equation}
\end{lem}
	\noindent\emph{Proof.}
	Clearly, we have
\begin{align*}
	\frac{2}{n+2}\binom{n-1}{\ell-1}\binom{n+\ell+1}{\ell}
	&=\left(1-\frac{n}{n+2}\right)\binom{n-1}{\ell-1}\binom{n+\ell+1}{\ell}\notag\\
	&=\binom{n-1}{\ell-1}\binom{n+\ell+1}{\ell}-\binom{n}{\ell}\binom{n+\ell+1}{n+2}.\notag
\end{align*}
This ends the proof.
\qed
	
Let $(x)_{0}=1$ and $(x)_{n}=x(x+1)\cdots(x+n-1)$ for all $n\geq 1$. In 2012, Guo and Zeng\cite{GuoZeng2012a} proved that for any given $k,\ell,a\in\bN$, there exist integers  $C_{0}{(\ell,a)},C_{1}{(\ell,a)}, \cdots,C_{a}{(\ell,a)}$ independent of $k$, such that
\begin{equation}\label{eq:k(k+1)}
	k^{a}(k+1)^{a}
	=\sum_{u=0}^{a}C_{u}{(\ell,a)}\prod_{v=1}^{u}\left(k(k+1)-(\ell+v-1)(\ell+v)\right)
\end{equation}
and then
\begin{equation*}
	k^{a}(k+1)^{a}\binom{k+\ell}{2\ell}
	=\sum_{u=0}^{a}C_{u}{(\ell,a)}\binom{
		k+\ell+u}{2\ell+2u}(2\ell+1)_{2u}.
\end{equation*}
Multiplying both sides by $\binom{2\ell}{\ell}$, we obtain
\begin{equation}\label{eq:(l+1)u}
	k^{a}(k+1)^{a}\binom{k+\ell}{2\ell}\binom{2\ell}{\ell}
	=\sum_{u=0}^{a}K_{u}{(\ell,a)}\binom{
		k+\ell+u}{2\ell+2u},
\end{equation}
where
\begin{equation}\label{eq:K_u}
	K_{u}{(\ell,a)}=C_{u}{(\ell,a)}\binom{2\ell+2u}{\ell+u}(\ell+1)_{u}^{2}.
\end{equation}
Equality \eqref{eq:k(k+1)} also shows
\begin{equation}\label{eq:C_{a}{(l,a)}}
	C_{a}{(\ell,a)}=1 \text{ for } \ a\geq0,
\end{equation}
\begin{equation}\label{eq:C_{0}{(l,a)}}
C_{0}{(\ell,a)}=\ell(\ell+1)\sum_{u=1}^{a}(-1)^{u-1}C_{u}{(\ell,a)}
\prod_{v=2}^{u}(\ell+v-1)(\ell+v)\quad \text{ for } \ a\geq1
\end{equation}
and
\begin{equation*}
	C_{1}{(\ell,a)}=\sum_{u=2}^{a}(-1)^{u}C_{u}{(\ell,a)}\sum_{j=1}^{u}
	\prod_{\substack{v\neq j \\ 1\leq v \leq u}}(\ell+v-1)(\ell+v)\quad \text{ for } \ a>1.
\end{equation*}
Thus one can see
\begin{lem}\label{lem:zhengchul}
Given $\ell\in\bZ^{+}$,  then
$\frac{C_{0}{(\ell,1)}}{\ell(\ell+1)}=1$.
When $a>1$, $\frac{C_{0}{(\ell,a)}}{\ell(\ell+1)}$ and $C_1(\ell,a)$ are even integers.
\end{lem}

In 2018, Chen and Guo \cite{chenguo2018} proved that for $k,i,j\in\bN$ and $h\geq1$, there exist integers $b_{i,t}^{(h)}(i\leq t \leq hi)$ divisible by  $\binom{t}{i}$ such that
\begin{equation}\label{eq:h-1}			\binom{k+i}{2i}^{h}\binom{2i}{i}=\sum_{t=i}^{hi}b_{i,t}^{(h)}\binom{k+t}{2t}\binom{2t}{t}
\end{equation}
and that 	
\begin{equation}\label{eq:reduction2-1}
\binom{k+i}{2i}\binom{2i}{i}\binom{k+j}{2j}\binom{2j}{j}
=\sum_{\ell=0}^{i+j}\binom{i+j}{i}\binom{j}{i+j-\ell}
                 \binom{\ell}{j}\binom{k+\ell}{2\ell}\binom{2\ell}{\ell}.
\end{equation}

By \eqref{eq:reduction2-1} and the method of induction, it is straight forward to prove the following Lemma.
\begin{lem}\label{lem:{k+l}{2l}}
For ${i_{j}}\in\bN$, $j=1,2,\cdots,m$, there exist integers $B_{i_{1},\ldots,i_{m}}^{(\ell)}$ such that
\begin{equation}\label{eq:m reduce}		
\prod_{j=1}^{m}\binom{k+i_{j}}{2i_{j}}\binom{2i_{j}}{i_{j}}
=\sum_{\ell=0}^{i_{1}+\cdots+i_{m}}B_{i_{1},\ldots,i_{m}}^{(\ell)}
      \binom{k+\ell}{2\ell}\binom{2\ell}{\ell}.
		\end{equation}
\end{lem}

\begin{lem}\label{C even}
Let $n, M, \ell\in\bZ^{+}$ and $I,i_{j}\in\bN$ for $1\leq j\leq 2M$.
We have
\begin{equation}\label{eq:even}
\sum_{\substack{i_{1}+\cdots+i_{2M}=I \\ 0\leq i_{1},\ldots,i_{2M}\leq n}}B_{i_{1},\ldots,i_{2M}}^{(\ell)}\equiv0\pmod 2.
\end{equation}
\end{lem}
\pf
We proceed by induction on $M$.
When $M=1$, by \eqref{eq:reduction2-1} we have
\begin{equation}\label{eq:i_1,i_2}
B_{i_{1},i_{2}}^{(\ell)}=
\binom{i_{1}+i_{2}}{i_{1}}\binom{i_{2}}{i_{1}+i_{2}-\ell}\binom{\ell}{i_{2}}.
\end{equation}
Note that $B_{i,I-i}^{(\ell)}=B_{I-i,i}^{(\ell)}$ if $0\leq i \leq I$.
When $I=2J+1$ is odd, clearly
\[
\sum_{\substack {i_{1}+i_{2}=I\\0\leq i_1,i_2\leq n}}B_{i_{1},i_{2}}^{(\ell)}=2\sum_{i=\max\{0,2J+1-n\}}^{J}B_{i,2J+1-i}^{(\ell)}.
\]
When $I=2J$ is even, then
\[
\sum_{\substack {i_{1}+i_{2}=I\\0\leq i_1,i_2\leq n}}B_{i_{1},i_{2}}^{(\ell)}
=2\sum_{i=\max\{0,2J-n\}}^{J-1}B_{i,2J-i}^{(\ell)}+B_{J,J}^{(\ell)}.
\]
Thus \eqref{eq:even} holds when $M=1$ as $B_{J,J}^{(\ell)}$ is even by definition.

Next we will assume the statement is true for $M$ and then show it is also satisfied for $M+1$.
Let $i=i_{2M+1}$, $j=i_{2M+2}$ in \eqref{eq:reduction2-1}.
Then multiplied the obtained equation with \eqref{eq:m reduce} leads to
\[
B_{i_{1},\ldots,i_{2M+2}}^{(\ell)}
=\sum_{\ell_{1}=0}^{i_{1}+\cdots+i_{2M}}
          B_{i_{1},\ldots,i_{2M}}^{(\ell_{1})}
      \sum_{\ell_{2}=0}^{i_{2M+1}+i_{2M+2}}B_{i_{2M+1},i_{2M+2}}^{(\ell_{2})}
           B_{\ell_{1},\ell_{2}}^{(\ell)}.
\]
Then we have
\begin{align*}
	&
	\sum_{\substack {i_{1}+\cdots+i_{2M+2}=I\\0\leq i_1,\ldots,i_{2M+2}\leq n}}B_{i_{1},\ldots,i_{2M+2}}^{(\ell)}\notag\\
	=&\sum_{I_{1}=0}^{I}
	  \sum_{\ell_{1}=0}^{I_{1}}
      \sum_{\substack {i_{1}+\cdots+i_{2M}=I_{1}\\0\leq i_1,\ldots,i_{2M}\leq n}}B_{i_{1},\ldots,i_{2M}}^{(\ell_{1})}
      \sum_{\ell_{2}=0}^{I-I_{1}}B_{\ell_{1},\ell_{2}}^{(\ell)}
      \sum_{\substack {i_{2M+1}+i_{2M+2}=I-I_{1}\\ 0\leq i_{2M+1},i_{2M+2}\leq n}}B_{i_{2M+1},i_{2M+2}}^{(\ell_{2})}
.\notag
\end{align*}
When $\ell_2>0$, we know
$
\sum\limits_{\substack {i_{2M+1}+i_{2M+2}=I-I_{1}\\ 0\leq i_{2M+1},i_{2M+2}\leq n}}B_{i_{2M+1},i_{2M+2}}^{(\ell_{2})}
$
is even by the induction which leads to the conclusion.
When $\ell_2=0$, one can see $B_{\ell_{1},0}^{(\ell)}=\binom{0}{\ell_1-\ell}$ equals $0$ unless $\ell_1=\ell>0$.
Then the conclusion follows from the fact that
$
 \sum\limits_{\substack {i_{1}+\cdots+i_{2M}=I_{1}\\0\leq i_1,\ldots,i_{2M}\leq n}}B_{i_{1},\ldots,i_{2M}}^{(\ell_{1})}
$
is even for any $\ell_1\in\bZ^{+}$.
\qed

\emph{Proof of \eqref{eq:1} in Theorem \ref{TM:Dn2}.}	
By the definition of $D_k^{(h)}(x)$, one can check that
\begin{align*}
		&\sum_{k=1}^{n}(-1)^{k}k^{a}(k+1)^{a}(2k+1)D_{k}^{(h)}(x)^{m}\notag\\
		=&\sum_{k=1}^{n}(-1)^{k}k^{a}(k+1)^{a}(2k+1)\left(\sum_{i=0}^{k}\binom{k+i}{2i}^{h}\binom{2i}{i}^{h}{x}^{i}\right)^{m}\notag\\
		=&\sum_{k=1}^{n}(-1)^{k}k^{a}(k+1)^{a}(2k+1)\sum_{0\leq i_{1},\ldots,i_{m}\leq k}\prod_{j=1}^{m}
		\binom{k+i_{j}}{2i_{j}}^{h}\binom{2i_{j}}{i_{j}}^{h}{x}^{i_{j}}\notag\\
	    =&\sum_{0\leq i_{1},\ldots,i_{m}\leq n}\sum_{k=0}^{n}(-1)^{k}k^{a}(k+1)^{a}(2k+1)\prod_{j=1}^{m}
		\binom{k+i_{j}}{2i_{j}}^{h}\binom{2i_{j}}{i_{j}}^{h}{x}^{i_{j}}\notag\\
=&\sum_{0\leq i_{1},\ldots,i_{m}\leq n}
\left(\prod_{j=1}^{m}\binom{2i_{j}}{i_{j}}^{h-1}x^{i_{j}}\right)
  \sum_{k=0}^{n}(-1)^{k}k^{a}(k+1)^{a}(2k+1)\prod_{j=1}^{m}
  \sum_{t=i_{j}}^{hi_{j}}b_{i_{j},t}^{(h)}
                         \binom{k+t}{2t}\binom{2t}{t}\notag
\end{align*}
with the help of equation \eqref{eq:h-1}.
By Lemma \ref{lem:{k+l}{2l}}, there exist integers $\tilde{B}_{i_{1},\ldots ,i_{m}}^{(\ell,h)}$ such that
\begin{equation*}\label{eq:Subsitute}
\prod_{j=1}^{m}\sum_{t=i_{j}}^{hi_{j}}b_{i_{j},t}^{(h)}
                         \binom{k+t}{2t}\binom{2t}{t}
=\sum_{\ell=0}^{hI}
       \tilde{B}_{i_{1},\ldots ,i_{m}}^{(\ell,h)}\binom{k+\ell}{2\ell}\binom{2\ell}{\ell},
\end{equation*}
where
$\tilde{B}_{i_{1},\ldots,i_{m}}^{(\ell,h)}=
\sum\limits_{\substack{i_j\leq t_j \leq hi_j\\ 1\leq j\leq m}}\left( \prod\limits_{s=1}^{m}b_{i_s,t_s}^{(h)}\right)B_{t_{1},\ldots,t_{m}}^{(\ell)}$ and $I=\sum\limits_{j=1}^{m}i_j$.
Then we have
\begin{align}\label{eq:F}
&\sum_{k=1}^{n}(-1)^{k}k^{a}(k+1)^{a}(2k+1)D_{k}^{(h)}(x)^{m}\notag\\
=&\sum_{0\leq i_{1},\ldots,i_{m}\leq n}
\left(\prod_{j=1}^{m}\binom{2i_{j}}{i_{j}}^{h-1}x^{i_{j}}\right)
  \sum_{\ell=0}^{hI}\tilde{B}_{i_{1},\ldots ,i_{m}}^{(\ell,h)}
  \sum_{k=0}^{n}(-1)^{k}k^{a}(k+1)^{a}(2k+1)\binom{k+\ell}{2\ell}\binom{2\ell}{\ell}\notag\\
=&\sum_{0\leq i_{1},\ldots,i_{m}\leq n}
	\left(\prod_{j=1}^{m}\binom{2i_{j}}{i_{j}}^{h-1}x^{i_{j}}\right)
	\sum_{\ell=0}^{hI}\tilde{B}_{i_{1},\ldots ,i_{m}}^{(\ell,h)}
	\sum_{u=0}^{a}K_{u}{(\ell,a)}
	\sum_{k=0}^{n}(-1)^{k}(2k+1)\binom{k+\ell+u}{2\ell+2u}\notag\\
	=&\sum_{0\leq i_{1},\ldots,i_{m}\leq n}\left(\prod_{j=1}^{m}\binom{2i_{j}-1}{i_{j}}^{h-1}x^{i_{j}}\right)
     \sum_{\ell=0}^{hI}\tilde{B}_{i_{1},\ldots,i_{m}}^{(\ell,h)}
     \sum_{u=0}^{a} 2^{m(h-1)} F_u(\ell,a,n)
\end{align}
with the help of equality \eqref{eq:(l+1)u} and Lemma \ref{lem: D_{k}^{(h)}(x)}.
Here $K_{u}{(\ell,a)}$ is defined in \eqref{eq:K_u} and
\[F_u(\ell,a,n)=K_{u}{(\ell,a)}(-1)^{n}(n+1+\ell+u)\binom{n+\ell+u}{2\ell+2u}.\]
As $m\geq 1$ and $h>1$, we have
\[\frac{2^{m(h-1)}F_0(\ell,a,n)}{n(n+1)(n+2)}
=(-1)^n\frac{C_{0}{(\ell,a)}}{\ell}\frac{2^{m(h-1)}}{n+2}
\binom{n+\ell+1}{\ell} \binom{n-1}{\ell-1}\in\bZ
\]
with the help of  Lemma \ref{lem:{2}{n+2}} and Lemma \ref{lem:zhengchul}.
When $u>0$, one can check that
\[
\frac{F_u(\ell,a,n)}{n(n+1)(n+2)}=(-1)^n
C_{u}{(\ell,a)}\binom{n+1+\ell+u}{n+2}\binom{n-1}{n-\ell-u}(\ell+1)^{2}_{u-1}\in\bZ.
\]
This completes the proof.
\qed

\emph{Proof of \eqref{eq:3Dk} in Theorem \ref{TM:Dn2}}.
When $m$ is odd, \eqref{eq:3Dk} follows directly from \eqref{eq: (2,n){n(n+1)(n+2)}}.

Next, we assume $m=2M$ is even and $I=\sum_{j=1}^{m}i_j$.
Let $h=1$ in \eqref{eq:F}.
Notice that $\tilde{B}_{i_{1},\ldots ,i_{m}}^{(\ell,1)}=B_{i_{1},\ldots,i_{m}}^{(\ell)}$ given in Lemma \ref{lem:{k+l}{2l}}.
Then \eqref{eq:F} leads to
\begin{align*}\label{eq: {k+l}{2l}}
	&\sum_{k=1}^{n}(-1)^{k}k^{a}(k+1)^{a}(2k+1)D_{k}(x)^{m}\notag\\
\equiv&\sum_{0\leq i_{1},\ldots,i_{m}\leq n}x^{I}
     \sum_{\ell=0}^{I}B_{i_{1},\ldots,i_{m}}^{(\ell)} F_0(\ell,a,n)\\
\equiv&(-1)^n n(n+1)(n+2)\sum_{0\leq i_{1},\ldots,i_{m}\leq n}x^{I}
     \sum_{\ell=1}^{I}B_{i_{1},\ldots,i_{m}}^{(\ell)} \frac{C_{0}{(\ell,a)}}{\ell}\frac{1}{n+2}
\binom{n+\ell+1}{\ell} \binom{n-1}{\ell-1}\\
=&(-1)^{n}n(n+1)(n+2)\sum_{I=0}^{nm}x^{I}\sum_{\ell=1}^{I}\frac{C_{0}{(\ell,a)}}{\ell}
      \frac{\binom{n+\ell+1}{\ell}\binom{n-1}{\ell-1}}{n+2}
      \sum_{\substack{i_1+\cdots+i_m=I \\ 0\leq i_{1},\ldots,i_{m}\leq n}}B_{i_{1},\ldots,i_{m}}^{(\ell)}\notag\\
      \equiv&0\pmod {n(n+1)(n+2)}\notag
\end{align*}
with the help of Lemma \ref{lem:{2}{n+2}} and \ref{C even}.
Hence, we have
\begin{equation*}
	\frac{1}{n(n+1)(n+2)}\sum_{k=1}^{n}(-1)^{k}k^{a}(k+1)^{a}(2k+1)D_{k}(x)^{m}\in\bZ[x].
\end{equation*}The proof of \eqref{eq:3Dk} is now completed.
   \qed

\section{Congruences involving $S_{k}^{(h)}(x)^{m}$}
In this section, we will study the arithmetic properties of $S_{k}^{(h)}(x)^{m}$.
\begin{theo}\label{th:Sn}
	Let $S_k^{(h)}(x)$ be given as in \eqref{eq:S_{k}^{(h)}(x)}.
For any ${h,m,n,a}\in\bZ^{+}$, we have
\begin{equation}\label{eq: new(2,n)}
\frac{(2,n)}{n(n+1)(n+2)}
\sum_{k=1}^{n}k^{a}(k+1)^{a}(2k+1)S_{k}^{(h)}(x)^{m}\in\bZ[x],
\end{equation}
and
\begin{equation}\label{eq:new(2,m-1,n)} \frac{(2,m-1,n)}{n(n+1)(n+2)}\sum_{k=1}^{n}(-1)^{k}k^{a}(k+1)^{a}(2k+1)S_{k}^{(h)}(x)^{m}\in\bZ[x].
\end{equation}
\end{theo}
To prove Theorem \ref{th:Sn}, we first establish some lemmas needed.
\begin{lem}\label{lem:S{k}^{(h)}(x)}
	For $n,u,\ell\in\bN$, we have
	\begin{equation}\label{eq:k+l+u}
		\sum_{k=0}^{n}(2k+1)\binom{k+\ell+u}{2\ell+2u}
		=\frac{(n+1)(n+1+\ell+u)\binom{n+\ell+u}{2\ell+2u}}{\ell+1+u}.
	\end{equation}
\end{lem}
\noindent\emph{Proof.}
By Gosper’s algorithm, we get
\begin{equation}\label{eq:-(k-l)(k^{2}-l-1)}
	(2k+1)\binom{k+\ell+u}{2\ell+2u}
=\Delta_{k}\left(\frac{k(k-\ell-u)\binom{k+\ell+u}{2\ell+2u}}{\ell+1+u}\right).
\end{equation}
Summing both sides of equality \eqref{eq:-(k-l)(k^{2}-l-1)} over
$k$ from $0$ to $n$ leads to the equality \eqref{eq:k+l+u}.
\qed

In Section 2, we find that \eqref{eq:reduction2-1} plays an crucial role in the proofs.
The following Lemma provides a similar reduction on products of $\frac{1}{i+1}\binom{k+i}{2i}$.
\begin{lem}\label{lem: 1/i+1}
  	For ${i,j,k}\in\bN$, we have
\begin{equation}\label{eq:{1}{i+1}}
   \frac{1}{i+1}\binom{k+i}{2i}\binom{2i}{i}\frac{1}{j+1}\binom{k+j}{2j}\binom{2j}{j}\\
  =\sum_{\ell=0}^{i+j}A_{i,j}^{(\ell)}
            \frac{1}{\ell+1}\binom{k+\ell}{2\ell}\binom{2\ell}{\ell},
\end{equation}
where
\begin{equation}\label{eq:A_{i,j}}
A_{i,j}^{(\ell)}=\frac{1}{j}\binom{i+j}{i+1}\binom{j}{i+j-\ell}\binom{\ell+1}{\ell-j}\in\bZ
\text{ for any } i,j,\ell\in\bN
\end{equation}
and $A_{i,j}^{(\ell)}=A_{j,i}^{(\ell)}$.
\end{lem}
  \noindent\emph{Proof.}
  A special case of the Pfaﬀ-Saalsch\"utz identity (see \cite[p.44, Exercise 2.d]{RP1997}) reads
  \begin{equation}\label{eq:PSidentity}
  	\binom{x+a}{b}\binom{y+b}{a}=\sum_{i=0}^{a}\binom{x+y+i}{i}\binom{y}{a-i}\binom{x}{b-i}.
  \end{equation}
  Taking $x=k$, $y=j$, $a=i$, $b=k-j$ and $i=\ell^{'}$ in \eqref{eq:PSidentity}, we get
  \begin{equation}\label{eq:(j+i)! l^{'}!}
  	\binom{k+i}{2i}\binom{2i}{i}=\sum_{\ell^{'}=0}^{i}\frac{(j+i)! \ell^{'}!}{(j+\ell^{'})! (i)!}\binom{j}{i-\ell^{'}}\binom{k-j}{\ell^{'}}\binom{k+j+\ell^{'}}{\ell^{'}}.
  \end{equation}
  Multiplying both sides of \eqref{eq:(j+i)! l^{'}!} by $\frac{1}{(i+1)(j+1)}\binom{k+j}{2j}\binom{2j}{j}$, we obtain
  \begin{align*}
  	&\frac{1}{i+1}\binom{k+i}{2i}\binom{2i}{i}\frac{1}{j+1}\binom{k+j}{2j}\binom{2j}{j}\notag\\
  	=&\sum_{\ell^{'}=0}^{i}\frac{(j+i)! \ell^{'}!}{(j+\ell^{'})! (i+1)!}\binom{j}{i-\ell^{'}}\binom{k-j}{\ell^{'}}\binom{k+j+\ell^{'}}{\ell^{'}}\frac{1}{j+1}\binom{k+j}{2j}\binom{2j}{j}\notag\\
  	=&\sum_{\ell^{'}=0}^{i}\binom{i+j}{i+1}\binom{j}{i-\ell^{'}}\binom{j+\ell^{'}+1}{\ell^{'}}\frac{1}{j(j+\ell^{'}+1)}\binom{k+j+\ell^{'}}{2j+2\ell^{'}}\binom{2j+2\ell^{'}}{j+\ell^{'}}.
\end{align*}
Taking $\ell=j+\ell^{'}$, we get equality \eqref{eq:{1}{i+1}}.
The fact that $A_{i,j}^{(\ell)}$ is integral follows by the decomposition
$
A_{i,j}^{(\ell)}= \binom{i+j}{i+j-\ell}\binom{\ell}{j}\binom{\ell+1}{i+1}
    -\binom{i+j}{i+1}\binom{j}{i+j-\ell}\binom{\ell+1}{\ell-j}.
$
\qed

Utilizing \eqref{eq:{1}{i+1}} repeatedly gives the following Lemma.
\begin{lem}\label{lem:{1}{l+1}{k+l}{2l}}
For ${i_{j}}\in\bN$, $j=1,2,\cdots,m$, there exist integers $A_{i_{1},\ldots,i_{m}}^{(\ell)}$ such that
\begin{equation}\label{m reduce(l+1)}
\prod_{j=1}^{m}\frac{1}{i_{j}+1}\binom{k+i_{j}}{2i_{j}}\binom{2i_{j}}{i_{j}}	=\sum_{\ell=0}^{i_{1}+\cdots+i_{m}}A_{i_{1},\ldots,i_{m}}^{(\ell)}
\frac{1}{\ell+1}\binom{k+\ell}{2\ell}\binom{2\ell}{\ell}
	\end{equation}
and $A_{i'_{1},\ldots,i'_{m}}^{(\ell)}=A_{i_{1},\ldots,i_{m}}^{(\ell)}$ when
$i'_{1}\cdots i'_{m}$ is a permutation of $i_{1}\cdots i_{m}$.
\end{lem}
\emph{Proof of \eqref{eq: new(2,n)} in Theorem \ref{th:Sn}.}
By the definition of $S_k^{(h)}(x)$, one can check that
\begin{align*}
	&\sum_{k=1}^{n}\varepsilon^{k}k^{a}(k+1)^{a}(2k+1)S_{k}^{(h)}(x)^{m}\notag\\
	=&\sum_{k=1}^{n}\varepsilon^{k}k^{a}(k+1)^{a}(2k+1)\left(\sum_{i=0}^{k}\frac{1}{(i+1)^{h}}\binom{k+i}{2i}^{h}\binom{2i}{i}^{h}{x}^{i}\right)^{m}\notag\\
	=&\sum_{k=1}^{n}\varepsilon^{k}k^{a}(k+1)^{a}(2k+1)\sum_{0\leq i_{1},\ldots,i_{m}\leq k}\prod_{j=1}^{m}
	\frac{1}{(i_{j}+1)^{h}}\binom{k+i_{j}}{2i_{j}}^{h}\binom{2i_{j}}{i_{j}}^{h}{x}^{i_{j}}\notag\\
	=&\sum_{0\leq i_{1},\ldots,i_{m}\leq n}
	\left(\prod_{j=1}^{m}C_{i_{j}}^{h-1}{x}^{i_{j}}\right)
	\sum_{k=0}^{n}\varepsilon^{k}k^{a}(k+1)^{a}(2k+1)\prod_{j=1}^{m}
	\frac{1}{i_{j}+1}\binom{k+i_{j}}{2i_{j}}^{h}\binom{2i_{j}}{i_{j}}\notag\\
	=&\sum_{0\leq i_{1},\ldots,i_{m}\leq n}
	\left(\prod_{j=1}^{m}C_{i_{j}}^{h-1}{x}^{i_{j}}\right)
	\sum_{k=0}^{n}\varepsilon^{k}k^{a}(k+1)^{a}(2k+1)
	\prod_{j=1}^{m}
	\sum_{t=i_{j}}^{hi_{j}}\frac{b_{i_{j},t}^{(h)}}{i_{j}+1}\binom{k+t}{2t}\binom{2t}{t}.
\end{align*}
with the help of equality \eqref{eq:h-1}, where $\varepsilon=\pm1$.
Since $\binom{t}{i_{j}}\mid b_{i_{j},t}^{(h)}$.
Then
\[
a_{i_{j},t}^{(h)}
=b_{i_{j},t}^{(h)}\frac{t+1}{i_{j}+1}
=\frac{b_{i_{j},t}^{(h)}}{\binom{t}{i_{j}}}\binom{t+1}{i_{j}+1}
\]
are integers.
Thus
\begin{align*}
	\prod_{j=1}^{m}
	\sum_{t=i_{j}}^{hi_{j}}\frac{b_{i_{j},t}^{(h)}}{i_{j}+1}\binom{k+t}{2t}\binom{2t}{t}
	=\prod_{j=1}^{m}
	\sum_{t=i_{j}}^{hi_{j}}a_{i_{j},t}^{(h)}\binom{k+t}{2t}\binom{2t}{t}\frac{1}{t+1}.
\end{align*}
By Lemma \ref{lem:{1}{l+1}{k+l}{2l}}, one can see
\begin{equation}\label{m reduce(l+1),h}
	\prod_{j=1}^{m}\sum_{t=i_{j}}^{hi_{j}}a_{i_{j},t}^{(h)}
	\binom{k+t}{2t}\binom{2t}{t}\frac{1}{t+1}
	=\sum_{\ell=0}^{hI}
	\tilde{A}_{i_{1},\ldots ,i_{m}}^{(\ell,h)}\binom{k+\ell}{2\ell}\binom{2\ell}{\ell}\frac{1}{\ell+1},
\end{equation}
where $I=\sum\limits_{j=1}^{m}i_j$ and
\begin{equation}\label{eq:Atilde}
\tilde{A}_{i_{1},\ldots,i_{m}}^{(\ell,h)}=
\sum\limits_{\substack{i_j\leq t_j \leq hi_j\\1\leq j\leq m}}\left(\prod\limits_{s=1}^{m}a_{i_s,t_s}^{(h)}\right)A_{t_{1},\ldots,t_{m}}^{(\ell)}.
\end{equation}
Then, we get
\begin{align}\label{eq:k^{a}(k+1)^{a}(2k+1)S_{k}^{(h)}}
	&\sum_{k=1}^{n}\varepsilon^{k}k^{a}(k+1)^{a}(2k+1)S_{k}^{(h)}(x)^{m}\notag\\
	=&\sum_{0\leq i_{1},\ldots,i_{m}\leq n}
	\left(\prod_{j=1}^{m}C_{i_{j}}^{h-1}{x}^{i_{j}}\right)
	\sum_{\ell=0}^{hI}\tilde{A}_{i_{1},\ldots ,i_{m}}^{(\ell,h)} \sum_{k=0}^{n}\varepsilon^{k}k^{a}(k+1)^{a}(2k+1)\frac{\binom{k+\ell}{2\ell}\binom{2\ell}{\ell}}{\ell+1}\notag\\
	=&\sum_{0\leq i_{1},\ldots,i_{m}\leq n}
	\left(\prod_{j=1}^{m}C_{i_{j}}^{h-1}x^{i_{j}}\right)
	\sum_{\ell=0}^{hI}\tilde{A}_{i_{1},\ldots ,i_{m}}^{(\ell,h)}
	\sum_{u=0}^{a}\frac{K_{u}{(\ell,a)}}{\ell+1}
	\sum_{k=0}^{n}\varepsilon^{k}(2k+1)\binom{k+\ell+u}{2\ell+2u}
\end{align}
with the help of \eqref{eq:(l+1)u}.

Let $\varepsilon=1$.
Substituting equality  \eqref{eq:k+l+u} into \eqref{eq:k^{a}(k+1)^{a}(2k+1)S_{k}^{(h)}}, we have
\begin{align}\label{eq:Ku/l+1}
&\sum_{k=1}^{n}k^{a}(k+1)^{a}(2k+1)S_{k}^{(h)}(x)^{m}\notag\\
=&\sum_{0\leq i_{1},\ldots,i_{m}\leq n}\left(\prod_{j=1}^{m}C_{i_{j}}^{h-1}x^{i_{j}}\right)
  \sum_{\ell=0}^{hI}\tilde{A}_{i_{1},\ldots    ,i_{m}}^{(\ell,h)}
  \sum_{u=0}^{a}G_u^{(1)}(\ell,a,n),
\end{align}
where $G_u^{(1)}(\ell,a,n)=\frac{K_{u}{(\ell,a)}(n+1)(n+1+\ell+u)}{(\ell+1)(\ell+1+u)}
       \binom{n+\ell+u}{2\ell+2u}$.
It is easy to check that
\[\frac{G_0^{(1)}(\ell,a,n)}{n(n+1)}=
      \frac{C_{0}{(\ell,a)}}{\ell(\ell+1)}\binom{n+\ell+1}{\ell+1}\binom{n-1}{\ell-1}\in\bZ,
      \]
\[\frac{G_1^{(1)}(\ell,a,n)}{n(n+1)}=
      C_{1}{(\ell,a)}\binom{n+\ell+2}{\ell+2}\binom{n-1}{\ell}\in\bZ
      \]
and when $u\geq 2$,
\[\frac{G_u^{(1)}(\ell,a,n)}{n(n+1)}=
     C_{u}{(\ell,a)}\binom{n+1+\ell+u}{\ell+u+1}
                    \binom{n-1}{\ell+u-1}{(\ell+1)}_{u}{(\ell+2)}_{u-2}\in\bZ.
      \]
Thus, we obtain
\begin{equation}\label{eq:n(n+1) divides Sk}
	n(n+1)\mid \sum_{k=1}^{n}k^{a}(k+1)^{a}(2k+1)S_{k}^{(h)}(x)^{m}.
\end{equation}
Let $n$ be $n+1$ gets
\begin{equation}\label{eq:(n+1)(n+2) divides Sk}
	(n+1)(n+2)\mid \sum_{k=1}^{n}k^{a}(k+1)^{a}(2k+1)S_{k}^{(h)}(x)^{m}.
\end{equation}
Combining \eqref{eq:n(n+1) divides Sk} and \eqref{eq:(n+1)(n+2) divides Sk} gives
\eqref{eq: new(2,n)}.
This concludes the proof.
\qed

Next, we continue to establish the following lemmas to prove \eqref{eq:new(2,m-1,n)}.
In \cite{Sun2018}, Sun defined
\[
w(n,\ell)=\frac{1}{\ell}\binom{n-1}{\ell-1}\binom{n+\ell}{\ell-1}
	=\binom{n-1}{\ell-1}\binom{n+\ell}{\ell}-\binom{n}{\ell}\binom{n+\ell}{\ell-1}\in\bZ.
\]
\begin{lem}\label{lem:{l=0}^{2b}w(n,l+1)}
	Let ${n,b}\in\bN$, we obtain
	\begin{align}\label{eq:{l=0}^{2b}w(n,l+1)}
		\sum_{\ell=0}^{2b}w(n,\ell+1)\equiv1\pmod{2}
	\end{align}
and
\begin{equation}\label{eq:H(n,l)0}
\sum_{\ell=0}^{2b}H(n,\ell)\equiv 0\pmod{2}.
\end{equation}
Here
\begin{equation}\label{eq:H(n,l)dingyi}
	H(n,\ell)=w(n,\ell+1)+w(n+1,\ell+1).
\end{equation}
\end{lem}
\noindent\emph{Proof.}
We proceed by mathematical induction on $b$.
It is clearly that \eqref{eq:{l=0}^{2b}w(n,l+1)} holds for $b=0$ since $w(n,1)=1$. When $b\geq 1$, suppose that the statement is true for $b-1$,  that is
\begin{equation*}
	\sum_{\ell=0}^{2b-2}w(n,\ell+1)\equiv1\pmod{2}.
\end{equation*}
Next, we will show it is also satisfied for $b$. Notice that
\begin{align*}
	&w(n,2b)+w(n,2b+1)\notag\\ 	
	=&\frac{1}{2b}\binom{n-1}{2b-1}\binom{n+2b}{2b-1}+\frac{1}{2b+1}\binom{n-1}{2b}\binom{n+2b+1}{2b}\notag\\
	=&\binom{n-1}{2b}\binom{n+2b}{2b-1}\frac{n(n+1)}{2b(n-2b)(2b+1)}\notag\\
	=&\binom{n+2b}{n-2b}C_{2b}\notag\\
	\equiv&0\pmod{2}.\notag
\end{align*}
Hence, we get
\begin{equation*}
	\sum_{\ell=0}^{2b}w(n,\ell+1)=\sum_{\ell=0}^{2b-2}w(n,\ell+1)+\sum_{\ell=2b-1}^{2b}w(n,\ell+1)\equiv1\pmod{2}.
\end{equation*}
Then congruence \eqref{eq:H(n,l)0} follows directly by the definition. \qed

\begin{lem}\label{lem:2^a-1}
Let $A_{J,J}^{(\ell)}$ given by \eqref{eq:A_{i,j}}.
Then for $J\leq \ell\leq 2J$, we have $A_{J,J}^{(\ell)}$ is odd if and only if $J=2^c-1$ for some $c\in\bZ^{+}$.
\end{lem}
\noindent\emph{Proof.}
Clearly, we have
\begin{align}
A_{J,J}^{(\ell)}
&=\frac{1}{J}\binom{2J}{J+1}\binom{J}{2J-\ell}\binom{\ell+1}{\ell-J}\label{eq:form1}\\
&=\frac{2}{J+1}\binom{2J-1}{J}\binom{J}{2J-\ell}\binom{\ell+1}{\ell-J}\in\bZ
\label{eq:form2}.
\end{align}
	
For the sufficiency, let $J=2^c-1$, we can easily prove that
\begin{equation*}
    \binom{2J}{J+1}\equiv\binom{J}{2J-\ell}\equiv\binom{\ell+1}{\ell-J}\equiv1\pmod2
\end{equation*}
for any $J\leq\ell\leq2J$.
Then $A_{J,J}^{(\ell)}$ is odd by \eqref{eq:form1}.

For the necessity, suppose $A_{J,J}^{(\ell)}$ is odd,
it can be seen that $J$ is also odd by \eqref{eq:form2}.
Then we know $\binom{2J}{J+1}$ is odd by \eqref{eq:form1} which happens only when $J=2^c-1$.
This completes the proof.
\qed

\begin{lem}\label{lem:I odd}
Let $n,h,M\in\bZ^{+}$, $\ell, i_j\in\bN$ for $1\leq j\leq 2M$ and
$C_{i_j}=\frac{1}{i_j+1}\binom{2i_j}{i_j}$ is the Catalan number.
When $I$ is an odd integer, we have
\begin{equation*}
\sum_{\substack{i_{1}+\cdots+i_{2M}=I \\ 0\leq i_{1},\ldots,i_{2M}\leq n}}\tilde{A}_{i_{1},\ldots,i_{2M}}^{(\ell,h)}\prod_{j=1}^{2M}C_{i_{j}}^{h-1}\equiv0\pmod 2,
	\end{equation*}
where $\tilde{A}_{i_{1},\ldots,i_{2M}}^{(\ell,h)}$ is given by \eqref{eq:Atilde}.
\end{lem}
\pf
When $h>1$, we know there exists $i_j$ with $1\leq j\leq 2M$ is even since $I$ is odd and $i_{1}+\cdots+i_{2M}=I$. Then $C_{i_j}$ is even which concludes the proof.

In the following, we assume $h=1$ and proceed by mathematical induction on $M$.
Note that $\tilde{A}_{i_{1},\ldots,i_{2M}}^{(\ell,1)}=A^{(\ell)}_{i_{1},\ldots,i_{2M}}$.
When $M=1$, we need to prove
\begin{equation*}
	\sum_{\substack{i_{1}+i_{2}=I \\0\leq i_{1},i_{2}\leq n}} A_{i_{1},i_{2}}^{(\ell)}\equiv0 \pmod2,
\end{equation*}
which follows by the fact that $A_{i,I-i}^{(\ell)}=A_{I-i,i}^{(\ell)}$, $0\leq i \leq I$ and $I$ is odd.
When $M>1$, suppose that the statement is true for $M-1$, that is
\begin{equation*}
	\sum_{\substack{i_{1}+\cdots+i_{2M-2}=I \\ 0\leq i_{1},\ldots,i_{2M-2}\leq n}}A_{i_{1},\ldots,i_{2M-2}}^{(\ell)}\equiv0\pmod 2.
\end{equation*}
Next we will show it is also satisfied for $M$. Let $i=i_{2M-1}$, $j=i_{2M}$ in \eqref{eq:{1}{i+1}}.
 Then multiplied the obtained equation with \eqref{m reduce(l+1)} by taking $m=2M-2$ leads to
 \begin{equation*}
 	A_{i_{1},\ldots,i_{2M}}^{(\ell)}
 	=\sum_{\ell_{1}=0}^{i_{1}+\cdots+i_{2M-2}}A_{i_{1},\ldots,i_{2M-2}}^{(\ell_{1})}
 	\sum_{\ell_{2}=0}^{i_{2M-1}+i_{2M}}A_{i_{2M-1},i_{2M}}^{(\ell_{2})}
 	A_{\ell_{1},\ell_{2}}^{(\ell)}.
 \end{equation*}	
Then, we have
\begin{align*}
	&\sum_{\substack{i_{1}+\cdots+i_{2M}=I \\ 0\leq i_{1},\ldots,i_{2M}\leq n}}A_{i_{1},\ldots,i_{2M}}^{(\ell)}\notag\\
	=&
	\sum_{I_{1}=0}^{I}
	\sum_{\ell_{1}=0}^{I_{1}}
	\sum_{\ell_{2}=0}^{I-I_{1}}A_{\ell_{1},\ell_{2}}^{(\ell)} \sum_{\substack{i_{1}+\cdots+i_{2M-2}=I_{1} \\ 0\leq i_{1},\ldots,i_{2M-2}\leq n}}A_{i_{1},\ldots,i_{2M-2}}^{(\ell_{1})}
\sum_{\substack{i_{2M-1}+i_{2M}=I-I_1 \\ 0\leq i_{2M-1},i_{2M}\leq n}}A_{i_{2M-1},i_{2M}}^{(\ell_{2})}
.\notag
\end{align*}
Since $I$ is odd.
We know either $I_{1}$ or $I-I_1$ is odd.
Then $\sum\limits_{\substack{i_{1}+\cdots+i_{2M-2}=I_{1} \\ 0\leq i_{1},\ldots,i_{2M-2}\leq n}}A_{i_{1},\ldots,i_{2M-2}}^{(\ell_{1})}$ or $\sum\limits_{\substack{i_{2M-1}+i_{2M}=I-I_1 \\ 0\leq i_{2M-1},i_{2M}\leq n}}A_{i_{2M-1},i_{2M}}^{(\ell_{2})}$ is even.
In any cases, we have
\[
\sum\limits_{\substack{i_{1}+\cdots+i_{2M}=I \\ 0\leq i_{1},\ldots,i_{2M}\leq n}}A_{i_{1},\ldots,i_{2M}}^{(\ell)}\equiv0\pmod{2}.
\]
This ends the proof.
\qed

\begin{lem}\label{lem:I is even}
Given $M,n,h\in\bZ^{+}$, $j_s\in\bN$ for $1\leq s\leq M$.
Let $I,e$ be even integers,
$\tilde{A}_{j_{1},j_{2},\ldots,j_{M}}^{(\ell,h,2)}
=\tilde{A}_{j_{1},j_{2},\ldots,j_{M},j_{1},j_{2},\ldots,j_{M}}^{(\ell,h)}$ given by \eqref{eq:Atilde} for simplicity and $H(n,\ell)$  defined in \eqref{eq:H(n,l)dingyi}.
We have
	\begin{equation*}
		\sum_{\ell=0}^{e}H(n,\ell)
		\sum_{\substack{j_{1}+j_{2}+\cdots+j_{M}=\frac{I}{2} \\ 0\leq j_{1},j_{2},\ldots,j_{M}\leq n}}\tilde{A}_{j_{1},j_{2},\ldots,j_{M}}^{(\ell,h,2)}	\left(\prod_{s=1}^{M}C_{j_{s}}^{h-1}\right)^{2}\equiv0\pmod2.
	\end{equation*}
\end{lem}
\noindent\emph{Proof.}
By \eqref{m reduce(l+1),h}, we obtain
%
\begin{align*}
\tilde{A}_{j_{1},j_{2},\ldots,j_{M}}^{(\ell,h,2)}
&=\sum_{\ell_{1}=0}^{\frac{hI}{2}}\tilde{A}_{j_{1},\ldots,j_{M}}^{(\ell_{1},h)}
\sum_{\ell_{2}=0}^{\frac{hI}{2}}\tilde{A}_{j_{1},\ldots,j_{M}}^{(\ell_{2},h)}
A_{\ell_{1},\ell_{2}}^{(\ell)}\\
&\equiv\sum_{\ell_{3}=0}^{\frac{hI}{2}}
\left(\tilde{A}_{j_{1},\ldots,j_{M}}^{(\ell_{3},h)}\right)^{2}
A_{\ell_{3},\ell_{3}}^{(\ell)}\pmod2
\end{align*}
since $A^{(\ell)}_{\ell_{1},\ell_{2}}=A^{(\ell)}_{\ell_{2},\ell_{1}}$.
Then
\begin{align*}
	&\sum_{\ell=0}^{e}H(n,\ell)
	\sum_{\substack{j_{1}+j_{2}+\cdots+j_{M}=\frac{I}{2} \\ 0\leq j_{1},j_{2},\ldots,j_{M}\leq n}}\tilde{A}_{j_{1},j_{2},\ldots,j_{M}}^{(\ell,h,2)}\left(\prod_{s=1}^{M}C_{j_{s}}^{h-1}\right)^{2}\notag\\
	\equiv&\sum_{\ell=0}^{e}H(n,\ell)
	\sum_{\substack{j_{1}+j_{2}+\cdots+j_{M}=\frac{I}{2} \\ 0\leq j_{1},j_{2},\ldots,j_{M}\leq n}} \sum_{\ell_{3}=0}^{\frac{hI}{2}}
\left(\tilde{A}_{j_{1},\ldots,j_{M}}^{(\ell_{3},h)}\right)^{2}
A_{\ell_{3},\ell_{3}}^{(\ell)}\left(\prod_{s=1}^{M}C_{j_{s}}^{h-1}\right)^{2}\notag\\
	=&\sum_{\substack{j_{1}+j_{2}+\cdots+j_{M}=\frac{I}{2} \\ 0\leq j_{1},j_{2},\ldots,j_{M}\leq n}} \sum_{\ell_{3}=0}^{\frac{hI}{2}}
\left(\tilde{A}_{j_{1},\ldots,j_{M}}^{(\ell_{3},h)}\right)^{2}
\left(\prod_{s=1}^{M}C_{j_{s}}^{h-1}\right)^{2}
\sum_{\ell=\ell_{3}}^{\min\{e,2\ell_{3}\}}H(n,\ell)A_{\ell_{3},\ell_{3}}^{(\ell)}\\
\equiv&\sum_{\substack{j_{1}+j_{2}+\cdots+j_{M}=\frac{I}{2} \\ 0\leq j_{1},j_{2},\ldots,j_{M}\leq n}}
\sum_{c=0}^{\log_2(\frac{hI}{2}+1)}
\left(\tilde{A}_{j_{1},\ldots,j_{M}}^{(2^c-1,h)}\right)^{2}
\left(\prod_{s=1}^{M}C_{j_{s}}^{h-1}\right)^{2}
	\sum_{\ell=2^c-1}^{\min\{e,2(2^c-1)\}}H(n,\ell)\notag\\
=&\sum_{\substack{j_{1}+j_{2}+\cdots+j_{M}=\frac{I}{2} \\
0\leq j_{1},j_{2},\ldots,j_{M}\leq n}}	\sum_{c=0}^{\log_2(\frac{hI}{2}+1)}\left(\tilde{A}_{j_{1},\ldots,j_{M}}^{(2^c-1,h)}
\right)^{2}\left(\prod_{s=1}^{M}C_{j_{s}}^{h-1}\right)^{2}\left( \sum_{\ell=0}^{\min\{e,2(2^c-1)\}}	H(n,\ell)-\sum_{\ell=0}^{2^c-2}	H(n,\ell)\right)\notag\\
\equiv&0\pmod2
\end{align*}
with the help of Lemma \ref{lem:2^a-1} and Lemma \ref{lem:{l=0}^{2b}w(n,l+1)}.
This concludes the proof.
\qed

\emph{Proof of \eqref{eq:new(2,m-1,n)} in Theorem \ref{th:Sn}.}
When $m$ is odd.
Let $\varepsilon=-1$.
Substituting equality  \eqref{eq:2k+1} into \eqref{eq:k^{a}(k+1)^{a}(2k+1)S_{k}^{(h)}}, we have
\begin{align}\label{eq:G_u(-1)}
	&\sum_{k=1}^{n}(-1)^{k}k^{a}(k+1)^{a}(2k+1)S_{k}^{(h)}(x)^{m}\notag\\
	=&(-1)^n\sum_{0\leq i_{1},\ldots,i_{m}\leq n}\left(\prod_{j=1}^{m}C_{i_{j}}^{h-1}x^{i_{j}}\right)
	\sum_{\ell=0}^{hI}\tilde{A}_{i_{1},\ldots    ,i_{m}}^{(\ell,h)}
	\sum_{u=0}^{a}G_u^{(-1)}(\ell,a,n),
\end{align}
where $G_u^{(-1)}(\ell,a,n)=\frac{K_{u}{(\ell,a)}}{\ell+1}(n+1+\ell+u)\binom{n+\ell+u}{2\ell+2u}.$
Note that
\begin{equation}\label{eq:G_0^{(-1)}}
\frac{G_0^{(-1)}(\ell,a,n)}{n(n+1)}
=\frac{C_{0}{(\ell,a)}}{\ell(\ell+1)}\binom{n+\ell+1}{\ell}\binom{n-1}{\ell-1}\in\bZ,
\end{equation}
\begin{equation}\label{eq:G_1^{(-1)}}
\frac{G_1^{(-1)}(\ell,a,n)}{n(n+1)}
=C_{1}{(\ell,a)}\binom{n+\ell+2}{\ell+1}\binom{n-1}{\ell}\in\bZ,
\end{equation}
and when $u\geq 2$,
\begin{equation}\label{eq:G_u^{(-1)}}
\frac{G_u^{(-1)}(\ell,a,n)}{n(n+1)(n+2)}=
C_{u}{(\ell,a)}\binom{n+1+\ell+u}{\ell+u-1}
\binom{n-1}{\ell+u-1}{(\ell+1)_{u-1}}{(\ell+2)_{u-2}}\in\bZ.
\end{equation}
Thus, we obtain
\begin{equation}\label{eq:n(n+1) divides -Sk}
	n(n+1)\mid \sum_{k=1}^{n}(-1)^{k}k^{a}(k+1)^{a}(2k+1)S_{k}^{(h)}(x)^{m}.
\end{equation}
Let $n$ be $n+1$ gets
\begin{equation}\label{eq:(n+1)(n+2) divides -Sk}
	(n+1)(n+2)\mid \sum_{k=1}^{n}(-1)^{k}k^{a}(k+1)^{a}(2k+1)S_{k}^{(h)}(x)^{m}.
\end{equation}
The equations \eqref{eq:n(n+1) divides -Sk} and \eqref{eq:(n+1)(n+2) divides -Sk} show
\begin{equation*}
	\frac{(2,n)}{n(n+1)(n+2)}
	\sum_{k=1}^{n}(-1)^{k}k^{a}(k+1)^{a}(2k+1)S_{k}^{(h)}(x)^{m}\in\bZ[x].
\end{equation*}

Next we assume $m=2M$ is even.
If $a>1$, \eqref{eq:new(2,m-1,n)} follows from \eqref{eq:G_0^{(-1)}}--\eqref{eq:G_u^{(-1)}} together with Lemma \ref{lem:{2}{n+2}} and Lemma \ref{lem:zhengchul}.
If $a=1$, equalities \eqref{eq:C_{a}{(l,a)}} and \eqref{eq:C_{0}{(l,a)}} lead to
$C_1(\ell,1)=1$, $C_0(\ell,1)=\ell(\ell+1)$ and then
\[
\frac{G_0^{(-1)}(\ell,1,n)+G_1^{(-1)}(\ell,1,n)}{n(n+1)(n+2)}=\frac{H(n,\ell)}{2}.
\]
Thus quality \eqref{eq:G_u(-1)} derives
\begin{align}\label{eq:I=even,a=1}	
	&\frac{1}{n(n+1)(n+2)}\sum_{k=1}^{n}(-1)^{k}k(k+1)(2k+1)S_{k}^{(h)}(x)^{m}\notag\\
	=&\frac{(-1)^n}{2}\sum_{0\leq i_{1},\ldots,i_{m}\leq n}\left(\prod_{j=1}^{m}C_{i_{j}}^{h-1}x^{i_{j}}\right)
	\sum_{\ell=0}^{hI}\tilde{A}_{i_{1},\ldots ,i_{m}}^{(\ell,h)}H(n,\ell)\notag\\
	=&\frac{(-1)^n}{2}
	\sum_{I=0}^{nm}x^{I}
	\sum_{\ell=0}^{hI}H(n,\ell)
	\sum_{\substack{i_1+\cdots+i_m=I \\ 0\leq i_{1},\ldots,i_{m}\leq n}}\tilde{A}_{i_{1},\ldots ,i_{m}}^{(\ell,h)}\prod_{j=1}^{m}C_{i_{j}}^{h-1}.
\end{align}
If $I$ is odd, \eqref{eq:new(2,m-1,n)} holds with the help of Lemma \ref{lem:I odd}.
If $I$ is even, note that  $\binom{2M}{2d+1}$ is even for any $d\in\bN$ and $\tilde{A}_{i_{1},i_{2},\ldots,i_{2M}}^{(\ell,h)}
=\tilde{A}_{i_{1}',i_{2}',\ldots,i_{2M}'}^{(\ell,h)}$ when
$i'_{1}\cdots i'_{2M}$ is a permutation of $i_{1}\cdots i_{2M}$.
We obtain
\begin{equation}\label{eq:symmetry}
	\sum_{\substack{i_{1}+\cdots+i_{2M}=I \\ 0\leq i_{1},\ldots,i_{2M}\leq n}}\tilde{A}_{i_{1},\ldots,i_{2M}}^{(\ell,h)}\prod_{j=1}^{2M}C_{i_{j}}^{h-1}
	\equiv\sum_{\substack{j_{1}+j_{2}+\cdots+j_{M}=\frac{I}{2} \\ 0\leq j_{1},j_{2},\ldots,j_{M}\leq n}}\tilde{A}_{j_{1},j_{2},\ldots,j_{M}}^{(\ell,h,2)}\left(\prod_{s=1}^{M}C_{j_{s}}^{h-1}\right)^{2}\pmod2.
\end{equation}
Substituting \eqref{eq:symmetry} into \eqref{eq:I=even,a=1} concludes the proof with the help of Lemma \ref{lem:I is even}.
\qed

\end{document}